\documentclass{elsarticle}

\newcommand{\mn}{\mathbb{N}} %math N%
\newcommand{\var}[1]{\mathrm{var}\left(#1\right)}
\newcommand{\Esym}{\mathrm{E}}
\newcommand{\E}[1]{\Esym\left[#1\right]}
\newcommand{\Probsym}{\mathbb{P}}
\newcommand{\Prob}[1]{\Probsym\left[#1\right]}
\newcommand{\bQ}{\mathbf{Q}}
\newcommand{\bP}{\mathbf{P}}
\newcommand{\bD}{\mathbf{D}}
\newcommand{\bB}{\mathbf{B}}
\newcommand{\ident}{\mathbf{I}}
\newcommand{\bzero}{\mathbf{0}}

\newcommand{\bA}{\mathbf{A}}

\usepackage{graphicx,amssymb,amsmath}

\begin{document}

\begin{frontmatter}

\title{A discrete-time Markov modulated queuing system with batched
arrivals}
\author{Richard G. Clegg}
\ead{richard@richardclegg.org}
\address{Department of Electrical and Electronic Engineering, University
College, London WC1E 7JE}

\begin{abstract}
This paper examines a discrete-time queuing system with 
applications to telecommunications traffic.  The arrival process is
a particular Markov modulated process which belongs to the
class of discrete batched Markovian arrival processes.
The server process is a single server deterministic queue.  
A closed form exact
solution is given for the expected queue length and delay.  A
simple system of equations is given for the probability of the queue
exceeding a given length.  
\end{abstract}

\end{frontmatter}

% main text
\section{Introduction}

This paper provides a solution for the expected queue length and 
probability of a given queue length for a simple discrete-time queuing system.  
The queuing system in question processes one unit of work in one unit of 
time.  Work
arrives in integer units according to an arrival process with the 
following properties.
\begin{itemize}
\item The system has two states, {\em on\/} and {\em off\/}.  In an 
{\em off\/} state, no
work will arrive.
\item If the system is in an {\em off\/} state then, with probability $f_0$,
in the next time unit the system is also an {\em off\/} state. 
\item If the system is in an {\em off\/} state then, with probability $f_i$, 
in the next time unit the system will move to an {\em on\/} state which
will last for exactly $i$ time units and
then move to an {\em off\/} state. 
\item If the system is {\em on\/} then a non-zero integer number of work units 
will arrive in this time period.  The number of units of work which arrive
is an iid random variable with 
$g_n$ as the probability that exactly $n$ units of work will arrive
in this time period.
\end{itemize}  

This process can be modelled as a Markov-modulated process (MMP) which is
completely characterised by the parameters $f_i$ and $g_i$.  
This model has been studied by authors in several 
different areas, for example statistical 
physics \cite{wang1989}, the study of dynamical systems 
\cite{coelho2000} and modelling telecommunications traffic 
\cite{clegg2005},\cite{woolf2002}.  
In the last two papers, the model is considered as the
source of input to a network and hence it is natural to consider the
queuing properties of such a model.  In this paper expressions are
derived for the expected queue length at equilibrium and the 
probability
that the queue has a given length at equilibrium under certain 
natural restrictions (for example, that the utilisation of the system is less 
than one).  It is shown that 
the expected queue length is a function of only four variables,
the first and second moments of the parameters $f_i$ and $g_i$.  
From Little's law \cite{little1961}, the expected delay is proportional 
to the expected queue length.

In section \ref{sec:model} the model is introduced formally and some 
basic properties are derived.  The model is related to existing work
in queuing theory.  In section \ref{sec:solve} the model is 
solved to get equations for the expected queue length and the probability 
of a given
queue length.  In section
\ref{sec:simulation} the derived equations are compared with computer
simulations. 

\subsection*{Acknowledgements}

I would like to thank Dr Simon Eveson and Prof Maurice 
Dodson for their invaluable help and, in particular Dr Eveson for his
insight regarding the Sherman--Morrison formula.  This work was
partly conducted under the EPSRC grant GR/T10503.  Thanks are due
to the reviewers who helpfully referenced a number of related
papers of which I was unaware.

\section{A Markov queuing model}
\label{sec:model}

Consider the discrete time process described in the introduction.
The motivation behind this process
is the idea that the lengths of {\em off\/} periods are 
memoryless (the probability that
an {\em on\/} period begins is independent of the length of the 
{\em off\/} 
period so far) 
but {\em on\/} periods have lengths which are iid with an arbitrary
distribution (within certain feasible constraints).  
This arrival process is then used as the input
to a queue which can process one unit of work (either queued 
or newly arriving) per time period.

\begin{figure}[ht]
\includegraphics[width=12cm]{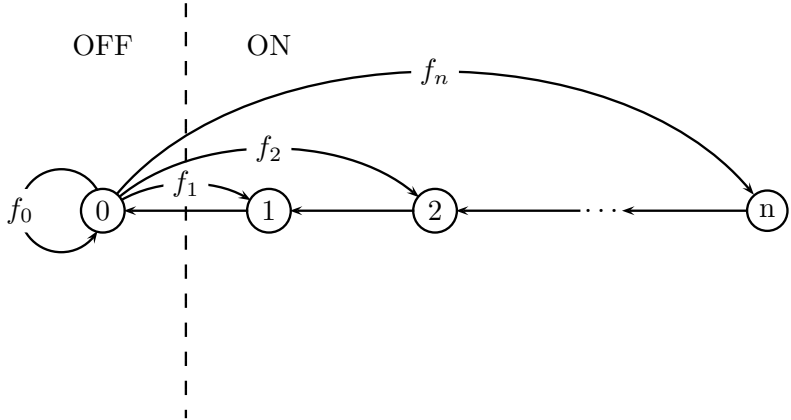}
\caption{The Markov chain for an MMP traffic model.}
\label{fig:cleggtop}
\end{figure}

The process $Y_t$, the arrival process to the system, can be modelled 
by an MMP.  The process has the underlying
Markov chain shown in figure \ref{fig:cleggtop}.  At any
time the chain is in an {\em on\/} state (all of the non-zero states)
then a non-zero number of arrivals will occur with a given 
probability.

The $n+1$ state Markov chain from figure \ref{fig:cleggtop} has
the transition matrix
\begin{equation*}
\bP= \left[
\begin{array}{cccccc}
f_0 & f_1 & \dots & f_{n-1} & f_n \\
1 & 0 & \dots & 0 & 0 \\
0 & 1 & \dots & 0 & 0 \\
\vdots & \vdots & \ddots & \vdots & \vdots \\
0 & 0  & \dots & 1 & 0 \\
\end{array}
\right].
\end{equation*}

Let $\{X_t: t \in \mn\}$ be the states
of a discrete-time, homogeneous Markov chain at time $t$.  This chain 
is shown in figure \ref{fig:cleggtop} and the parameters $f_i$ define
the transition probabilities of the chain.  It can be seen from the
structure of the chain that $f_i$ is the probability that, given the 
chain is currently in state zero, the next state will be $i$.
This could also be thought as the probability that, given the
chain is currently in an {\em off\/} state, the next state will begin
an {\em on\/} period of length exactly $i$ ($f_0$ is the probability the
next state will be an {\em off\/} state again).

Let $g_i: i \in \mn$ be the 
probability that exactly $i$ units of work arrive when the chain is
in an {\em on\/} (non-zero) state ($X_t > 0$).    
Let $\{Y_t: t \in \mn\}$ be a series derived from 
$X_t$ by the rule, that if $X_t = 0$ then $Y_t = 0$ and if $X_t > 0$ 
then
for all $i \in \{1,2,\ldots,m\}$ 
(where $m$ is the maximum possible value of $Y_t$) 
$\Prob{Y_t= i} = g_i$ (for an {\em on\/} state, 
the possibility that $Y_t = 0$ is excluded, that is, $g_0 = 0$).  

The server model used is a deterministic process where exactly one
unit of work is processed in one time unit.  

\subsection{Related work}

The arrival process belongs to the processes known as
discrete batch Markovian  arrival process (D-BMAP) and the 
queuing system is a subset of D-BMAP/D/1 queues.  
The BMAP itself was introduced by
Lucantoni \cite{lucantoni1991}.  Both the D-BMAP and the BMAP have 
previously been used as models of telecommunications traffic 
\cite{blondia1992,hermann1993,klemm2003}.
Details of the BMAP can be found in most modern books on queuing
theory \cite[Chapter 12]{breuer2005} 
and only a brief outline is given here.
It is also interesting to note that the arrival model considered here
is very similar to the batch renewal process studied by Fretwell
and Kouvatsos \cite{fretwell1996} in the context of
internet traffic.  In fact the system they study
is the one in figure \ref{fig:cleggtop} with on and off reversed.

The structure of a generic D-BMAP is as follows.  Let $\bP$ be
the transition matrix for a discrete time Markov chain with
state space $E= \mn \times {0, \ldots,n}$.   Each state
of the chain is a pair $(j,i)$ where $j$ is the {\em level\/} 
(the number of arrivals generated in that state) and $i$ is the 
{\em phase\/}.  The transition matrix has the structure
$$
\bP = \left[ 
\begin{array}{ccccc}
\bD_0 & \bD_1 & \bD_2 & \bD_3 & \cdots \\
\bzero & \bD_0 & \bD_1 & \bD_2  & \cdots \\
\bzero & \bzero & \bD_0 & \bD_1  & \cdots \\
\vdots & \vdots & \vdots & \vdots & \ddots
\end{array}
\right],
$$
where $\bD_i$ is, itself an $(n+1)\times(n+1)$ matrix.  The $(j,k)$th entry
in $\bD_i$ is the probability of a phase transition from phase $j$ to $k$
given level $i$.

The arrival process described in the previous section could be 
described as a D-BMAP where the phase is the state of the chain
in Figure \ref{fig:cleggtop} and the levels simply represent
the various {\em on\/} states.  The $\bD$ are
the $(n+1)\times(n+1)$ matrices
$$
\bD_0 = \left[
\begin{array}{cccc}
f_0 & 0 & \cdots & 0 \\
f_1 & 0 & \cdots & 0 \\
\vdots & \vdots & \ddots & \vdots \\
f_n & 0 & \cdots & 0
\end{array}
\right].
$$
and
$$
\bD_i = \left[
\begin{array}{ccccc}
0 & g_i & 0 & \cdots & 0 \\
0 & 0 & g_i & \cdots & 0 \\
\vdots & \vdots & \vdots & \ddots & \vdots \\
0 & 0 & 0 & \cdots & g_i \\
0 & 0 & 0 & \cdots & 0 \\
\end{array}
\right],
$$
for $i \in \{1,2,\ldots,m\}$.

Blondia and Cassals \cite{blondia1992} provide a method for solving the
D-BMAP/D/1/K queuing model which gives a solution in
terms of a recursive series of matrix equations.  The
complexity of calculation is given by the author as
$K^2M^3$ where $K$ is the buffer size (the possible
number of items in the queue) and $M$ is the 
number of phases ($n+1$ in the system described here).

In contrast the system studied here only works for
infinite buffers and gives an answer for the expected
queue length in closed form.  It gives an
answer for the probability of
the queue having a given length as a recursive system
of equations each with $n$ terms.  

\subsection{Basic properties of the system}

It is useful to define 
$\overline{f}= \sum_{i=1}^n if_i$ and $\overline{f^2}= 
\sum_{i=1}^n i^2 f_i$.
Similarly, define
$\overline{g} = \sum_{i=1}^m i g_i$ and also 
$\overline{g^2} = \sum_{i=1}^m i^2 g_i$.
Since $\overline{g}$ is the mean number of arrivals in an {\em on\/} 
state  then this
should be finite (otherwise the system will have an infinite mean 
arrival
rate).  This will be the case in all systems with $m$ finite.

Let $\pi_i$ be the equilibrium probability of state $i$.  This exists 
when the chain is ergodic.  It can be
easily shown that the finite chain is ergodic if $f_0 \in (0,1)$ and,
if this is the case, the equilibrium probability of state $i$, $\pi_i$ is 
given
by
\begin{equation*}
\pi_i =  \pi_0 \sum_{j=i}^n f_j.
\end{equation*}
and
\begin{equation*}
\pi_0 = 1 - \sum_{i=1}^n \pi_i = \frac{1}{1 + \overline{f}}
\end{equation*}
and rearranging gives
\begin{equation*}
\overline{f}= \frac{1 - \pi_0}{\pi_0}.
\end{equation*}

It is also of interest to consider the version of this chain where 
$n \rightarrow \infty$.  That is for any $N \in \mn$ there exists 
$i > N$ such that $f_i > 0$. 
The mean return time to the state zero is given by $1 + \overline{f}$ 
which
is finite if $\overline{f}$ is finite.  Therefore the infinite
chain is ergodic if $f_0 > 0$ and $\overline{f}$ is finite.  The 
infinite chain
is useful in studying long-range dependence \cite{barenco2004,clegg2005}.  
Obviously for finite chains $\overline{f}$ is finite.  

The output $Y_t$ is then used as input to a queuing system.  Let 
$Q_t$ represent
the number queuing at time $t$ and $Y_t$ represent the number of 
arrivals to the
system during $[t, t+1)$. One queued item is removed from the queue 
every time
unit if there is any work in the queue or if any arrives (if $Q_t + 
Y_t > 0$).
Therefore the queuing system is as follows
\begin{equation}
Q_{t+1} = \left[ Q_t + Y_t - 1 \right]^+,
\label{eqn:queuesys}
\end{equation}
where $[X]^+$ means $\max(0,X)$.

The quantity $1 - \pi_0 = \overline{f}/(1+\overline{f})$ 
represents the proportion of the time the system is in an {\em on }
state and
therefore $\overline{g}(1 - \pi_0)$ 
is the mean arrival rate $\lambda$.  Since the
system can output one unit of work per time period the queue
utilisation (proportion of time the queue is non-empty)
$\rho = \lambda$.  Both are given by,
\begin{equation}
\rho = \lambda = \frac{\overline{g}\overline{f}}{1+\overline{f}}.
\label{eqn:rho}
\end{equation}
In order that the queue does not grow forever it is obviously 
a necessary condition that $\rho < 1$.

To summarise, the model is specified by the parameters $f_i$ and $g_i$.
The requirements on these parameters are that $f_0 > 0$ 
(which guarantees the underlying MC is aperiodic), that 
$\overline{f}$ is finite and that
$\overline{g}\overline{f}/(1+\overline{f}) < 1$ (which in
turn is a requirement that $\overline{g}$ is finite).  
The first two requirements
together guarantee that the underlying MC is ergodic the third
requirement ensures that the queue utilisation is less than one.

\subsection{Notation}

The following notation is used in this paper and is gathered
here for convenience.
\begin{itemize}
\item $Y_t$ --- the number of arrivals to the system during $[t,t+1)$.
\item $X_t$ --- the state of the underlying Markov chain during 
$[t,t+1)$.
\item $Q_t$ --- the number queuing at time $t$.
\item $X, Y, Q$ --- the above quantities as random variables
at some time when the system is in equilibrium.
\item $\pi_i$ --- the equilibrium probability of the $i$th state of 
the chain.
\item $\bP$ --- the transition matrix for $X_t$.
\item $f_i$ --- the transition probabilities in $\bP$.
\item $n$ --- the highest possible value of $X_t$ (the highest 
numbered state in
$\bP$).
\item $\overline{f}, \overline{f^2}$ --- the first and second moments 
of $f_i$, 
$\sum_{i=1}^n i f_i$ and $\sum_{i=1}^n  i^2 f_i$.
\item $Q_i(z)$ --- the generating function $\sum_{q=0}^\infty \Prob{Q 
= q , X = i} z^q$.
\item $\bQ(z)$ --- the $n+1$ column vector 
$\left[Q_0(z), Q_1(z), \dots, Q_n(z)\right]^T$. 
\item $g_i$ --- the probability that an amount of work $i$ arrives in 
the
next time unit if the system is on, $\Prob{Y_t = i|X_t > 0}$.
\item $m$ --- the maximum possible value of $Y_t$ (the largest number 
of units
of work which may arrive in unit time).
\item $\overline{g}, \overline{g^2}$  --- the first and second 
moments of $g_i$, 
$\sum_{i=1}^m i g_i$ and $\sum_{i=1}^m  i^2 g_i$.
\item $g(z)$ --- the generating function $\sum_{i=1}^m g_i z^i$.
\item $b_0$ --- the boundary condition $\Prob{Q = 0 | X = 0}$.
\item $\rho$ --- the queue utilisation.
\item $\lambda$ --- the mean arrival rate.
\end{itemize}

\section{Solving the queuing model}
\label{sec:solve}

Let $\bQ(z) = \left[Q_0(z), Q_1(z), \dots, Q_n(z)\right]^T$ 
be a column vector of the generating function for the queue in 
each state of the chain.  That is, 
$$Q_i(z) = \sum_{q=0}^\infty \Prob{Q 
= q , X = i} z^q.$$
Consider a general MMP arrival process.
Let $A_i(z)$ be the generating function for the number of
arrivals if the underlying chain is in state $i$.  Let
$\bB = [b_0, b_1, \dots, b_n]^T$ be the $(n+1)$ column vector 
of boundary conditions, $b_i = \Prob{Y=0, Q=0 | X = i}$.
Following Li \cite{li1993} it can be shown that the queuing system of 
equation
(\ref{eqn:queuesys}) using a general MMP with transition matrix
$\bP$ as input implies
\begin{equation}
\bQ(z)= (z-1)[z \ident - \bP^T \mathbf{G}(z))]^{-1}\bP^T\bB,
\label{eqn:queuegen}
\end{equation}
where $\mathbf{G}(z) = \mathrm{diag} (A_0(z), A_1(z), \dots, A_n(z))$. 

For the MMP in question $\bB$ and $\mathbf{G}(z)$ have much
simpler forms.  In the {\em off} state the generating function
for arrivals is simply 1 (no arrivals occur with probability one) 
and in the {\em on } state the generating function is $g(z)$.  Therefore,
$\mathbf{G}(z) = \mathrm{diag} (1, g(z), g(z), \dots , g(z))$.
Since in the {\em on\/} state, the system always has at least one
arrival and in the {\em off\/} state the system always has no 
arrivals then $\bB$ is the $(n+1)$ column vector, 
$\bB= [b_0, 0, 0, \dots, 0]^T$ where $b_0 = \Prob{Y=0, Q=0 | X=0} =
\Prob{Q=0 | X= 0}$.  It is these simplifications
which make this system soluble.

For  the specific MMP being
studied \eqref{eqn:queuegen} can be rewritten,
\begin{equation*}
\bQ(z)= (z-1)(\bA + \bA')^{-1}\bP^T\bB,
\end{equation*}
where $\bA$ is the $(n+1) \times (n+1)$ matrix
\begin{equation*}
\bA= \left[
\begin{array}{ccccc}
z & - g(z) & 0 & 0 & \dots \\
0 & z & -g(z) & 0 & \dots \\
0 & 0 & z  & -g(z) & \dots \\
0 & 0 & 0 & z  & \dots \\
\vdots & \vdots & \vdots & \vdots & \ddots
\end{array}
\right],
\end{equation*}
and $\bA'$ is the rank one $(n+1) \times (n+1)$ matrix which can be 
written as $\bA'=uv$ where $u=[f_0,f_1,f_2,\dots, f_n]^T$ and 
$v=[-1,0,0,\dots]$.

The matrix $\bA$ can be inverted giving the $(n+1) \times (n+1)$ 
matrix
\begin{equation*}
\bA^{-1}=\frac{1}{z} \left[
\begin{array}{ccccc}
(g(z)/z)^0 & (g(z)/z)^1 & (g(z)/z)^2 & (g(z)/z)^3 & \dots \\
0 & (g(z)/z)^0 & (g(z)/z)^1  & (g(z)/z)^2  & \dots \\
0 & 0 & (g(z)/z)^0   & (g(z)/z)^1  & \dots \\
0 & 0 & 0 & (g(z)/z)^0   & \dots \\
\vdots & \vdots & \vdots & \vdots & \ddots
\end{array}
\right].
\end{equation*}
Note that 
since $g_0 = 0$ then 
$g(z)/z) = \sum_{i=1}^m g_i z^i/z = 
\sum_{i=1}^{m} g_i z^{i-1}$.
If $z \in [0,1]$ then
$\sum_{i=1}^{m} g_i z^{i-1} \leq \sum_{i=1}^m g_i = 1$.
Therefore $g(z)/z \in [0,1]$ if $z \in [0,1]$
Hence $(g(z)/z)^n$ 
remains
bounded as $n \rightarrow \infty$.
The 
Sherman--Morrison formula (see, for example, \cite{bernstein2000}) states
that,
\begin{equation*}
(\bA+\bA')^{-1}= \bA^{-1} - \frac{\bA^{-1}uv\bA^{-1}}{1+v\bA^{-1}u}.
\end{equation*}
Now, 
\begin{equation*}
1+v\bA^{-1}u=  1 - 1/z\sum_{i=0}^n (g(z)/z)^i f_i.
\end{equation*}
Define 
\begin{equation*}
a_i(z) = \frac{1}{z}\sum_{j=i}^n f_j 
\left(\frac{g(z)}{z}\right)^{j-i}.
\end{equation*}
Then
\begin{equation*}
(\bA + \bA')^{-1}= \left[ \frac{[1 - a_0(z)]\ident -
\bA^{-1}uv}{1- a_0(z)}
\right] \bA^{-1},
\end{equation*}
whence,
\begin{equation*}
\bQ(z)= (z - 1) \frac{[1-a_0(z)]
\ident -
\bA^{-1}uv }{1 - a_0(z)}
\bA^{-1}\bP^T\bB.
\end{equation*}
Multiplying the matrices gives
\begin{equation*}
Q_i(z)=\frac{z-1}{1-a_0(z)}b_0a_i(z).
\end{equation*}
From the definition of $Q_i(z)$, 
$\E{z^Q} = \sum_{i=0}^n Q_i(z)$ and therefore
\begin{equation}
\E{z^Q} = b_0   \frac{(z-1)\sum_{i=0}^n f_i\sum_{j=0}^i(g(z)/z)^j}
{z - \sum_{i=0}^n f_i (g(z)/z)^i} =
\frac{b_0 N(z)}{D(z)},
\label{eqn:ezq}
\end{equation}
where $N(z)= (z-1)\sum_{i=0}^n f_i\sum_{j=0}^i(g(z)/z)^j$ and 
$D(z) = z - \sum_{i=0}^n f_i (g(z)/z)^i$.

\subsection{Calculating $b_0$}

Now it is necessary to calculate $b_0$.  Note
that $\lim_{z \rightarrow 1} z^Q = 1$, whence
\begin{equation*}
b_0 = \lim_{z \rightarrow 1} \frac{D(z)}{N(z)}.
\end{equation*}
Since, $\lim_{z \rightarrow 1} D(z) = \lim_{z \rightarrow 1} 
N(z) = 0$, by L'H\^{o}pital's rule,
\begin{equation*}
b_0 = \lim_{z \rightarrow 1} \frac{D'(z)}{N'(z)}.
\end{equation*}
But
\begin{equation*}
D'(z) = 1 - (g(z)/z)' \sum_{i=1}^n i f_i (g(z)/z) ^{i-1},
\end{equation*}
and 
\begin{equation*}
(g(z)/z)'= z^{-2}\sum_{i=1}^m \left( i g_i z^i - g_i z^i \right).
\end{equation*}
Hence $\lim_{z \rightarrow 1} (g(z)/z)' = \overline{g} - 1$, 
which implies,
\begin{equation}
\lim_{z \rightarrow 1} D'(z) = 1 - (\overline{g} - 1)\sum_{i=1}^n i 
f_i = 1 + (1 - \overline{g})\overline{f}.
\label{eqn:dz}
\end{equation}
Similarly,
\begin{equation*}
N'(z) = \sum_{i=0}^n f_i\sum_{j=0}^i(g(z)/z)^j + 
(z-1) \sum_{i=0}^n (g(z)/z)' f_i\sum_{j=0}^i j(g(z)/z)^{j-1}.
\end{equation*}
Providing the sum at the right hand side remains finite (which it will
for all finite $n$) then
the factor of $(z-1)$ will cancel this term as $z \rightarrow 1$.
This gives
\begin{equation*}
\lim_{z \rightarrow 1} N'(z) = \sum_{i=0}^n (i+1)f_i = 1 + 
\overline{f}.
\end{equation*}
Finally, therefore,
\begin{equation}
b_0 = 1 - \frac{\overline{g}\overline{f}}{1+\overline{f}} =
 1 - \rho.
\label{eqn:b0}
\end{equation}

\subsection{Queuing results}
\label{sec:eq}

The next stage is to get a function for the expectation of the queue 
size.
\begin{equation*}
\lim_{z \rightarrow 1} \frac{ d \E{z^Q}}{dz} = 
\lim_{z \rightarrow 1} \sum_{q=1}^\infty q \Prob{Q=q} z^{q-1} =
\sum_{q=1}^\infty q \Prob{Q=q} =\E{Q}.
\end{equation*}
Since $b_0$ is constant, from \eqref{eqn:ezq},
\begin{equation*}
\E{Q} = b_0\lim_{z \rightarrow 1}\frac{d (N(z)/D(z))}{dz} = 
b_0\lim_{z \rightarrow 1}\frac{N'(z)D(z) - N(z)D'(z)}{D(z)^2}.
\end{equation*}
Similarly, $\lim_{z \rightarrow 1} D(z) = 0$ and, from 
L'H\^{o}pital's rule, 
\begin{equation*}
\E{Q} =
b_0\lim_{z \rightarrow 1}\frac{N''(z)D(z) - N(z)D''(z)}{2D(z)D'(z)}
\end{equation*}
and, since $\lim_{z \rightarrow 1} N(z)/D(z) = 1/b_0$,
\begin{equation}
\E{Q}= b_0 \lim_{z \rightarrow 1}\frac{ N''(z) - D''(z)}{2 D'(z)}.
\label{eqn:expq}
\end{equation}
It is now necessary to find expressions for $N''(z)$ and $D''(z)$.  
Firstly,
\begin{eqnarray*}
N''(z)  & = [2(g(z)/z)' + (z-1)(g(z)/z)''] \sum_{i=1}^n f_i 
\sum_{j=1}^i j(g(z)/z)^{j-1} \\
 & + (z-1) (g(z)/z)'^2 \sum_{i=1}^n f_i \sum_{j=1}^i j(j-1) 
(g(z)/z)^{j-2}.
\end{eqnarray*}
Hence, if all the
sums remain finite (which they will if $n$ is finite),
\begin{equation}
\lim_{z \rightarrow 1} N''(z) =[2(\overline{g} - 1)]\sum_{i=1}^n 
i(i+1)f_i/2  
= (\overline{g} - 1) ( \overline{f^2} + \overline{f}).
\label{eqn:n2z}
\end{equation}
Similarly
\begin{align*}
D''(z)=  &-(g(z)/z)'' \sum_{i=1}^n i f_i (g(z)/z)^{i-1} \\ 
& \quad- (g(z)/z)'^2 \sum_{i=1}^n i (i-1)f_i(g(z)/z)^{i-2},
\end{align*}
and
\begin{equation*}
\lim_{z \rightarrow 1} (g(z)/z)'' = \lim_{z \rightarrow 1} 
z^{-3}\sum_{i=1}^m (i^2 - 3i +2)g_i z^i
= \overline{g^2} - 3\overline{g} + 2,
\end{equation*}
therefore,
\begin{equation}
\lim_{z \rightarrow 1} D''(z) = [\overline{g} + \overline{g}^2 - 
\overline{g^2} -1]\overline{f} 
-(\overline{g} - 1)^2 \overline{f^2}.
\label{eqn:d2z}
\end{equation}

Substituting \eqref{eqn:dz}, \eqref{eqn:b0}, \eqref{eqn:n2z} and 
\eqref{eqn:d2z} into \eqref{eqn:expq} gives
%\begin{equation}
%\E{Q}= \frac
%{ \overline{g}(\overline{g} - 1)
%\left[ \pi_0^2 \left(\sum_{i=1}^ni^2 f_i\right)  - (1-\pi_0)^2 \right]
%+ (\overline{g^2} - \overline{g}^2)(1-\pi_0)}
%{2[1 - \overline{g}(1 - \pi_0)]}.
%\end{equation}
\begin{equation}
\E{Q}= \frac
{
\overline{g}(\overline{g} - 1)
\left[ \overline{f^2} -  \overline{f}^2\right]  +
\overline{f}(1 + \overline{f})
\left[ \overline{g^2} - \overline{g}^2\right]
}
{2(1 + \overline{f})[1 + \overline{f} -  \overline{f}\overline{g}]}.
\label{eqn:eqfinal}
\end{equation}

Note that $\overline{g} \geq 1$, $\overline{f^2} \geq \overline{f}^2$ 
and
$\overline{g^2} \geq \overline{g}^2$ by their respective definitions
and  $1 + \overline{f} > \overline{f}\overline{g}$ for a system with
utilisation less than one by equation \eqref{eqn:rho}.  All bracketed 
terms in the numerator are therefore positive or zero and the 
denominator is 
strictly positive.  Note, however that $\overline{f^2}$ and
$\overline{g^2}$ are only guaranteed finite for systems with finite 
$n$ and $m$
respectively.  
That
is to say that some systems with a mean arrival rate less than one
and an utilisation less than one will still
have an no finite value for the expected queue length.  An example
of such a system are the systems with the $f_i$ parameters 
given in \cite{wang1989,clegg2005}
which both have $\overline{f^2}$ as a non-convergent series.
Such systems are of interested to those studying long-range dependence
and heavy-tailed distributions.

Interpreting $\overline{f^2} -  \overline{f}^2$ and
$\overline{g^2} -  \overline{g}^2$ as the variance of $f$ and $g$ 
respectively 
then \eqref{eqn:eqfinal} could also be written as
\begin{equation*}
\E{Q}= \frac
{
\overline{g}(\overline{g} - 1)\var{f} +
\overline{f}(1 + \overline{f})\var{g}
}
{2(1 + \overline{f})^2(1 -  \rho)}.
\end{equation*}
From \eqref{eqn:rho}, the expected delay is given from Little's law
\begin{equation*}
\E{T} = \frac{\E{Q}}{\lambda} =  \frac
{
\overline{g}(\overline{g} - 1)\var{f} +
\overline{f}(1 + \overline{f})\var{g}
}
{2 (1+ \overline{f})^2 \rho (1 - \rho)}.
\end{equation*}

An implication of these equations is that, assuming the mean traffic
level is fixed (that is $\overline{f}$ and $\overline{g}$ cannot be
changed) then the queueing delay of the system would be minimised if
the variance in the lengths of the {\em on } periods was minimised and
the variance of the amount of traffic arriving in an {\em on } period
was minimised.  This has an interesting analogy to the well-known
Pollaczek-Khinchin result \cite{khinchin1932} that for an M/G/1 queue 
the waiting time is proportional to the variance in the service time.

\subsection{Finding the queue distribution function}
\label{sec:pq}

In order to be able to ask questions about, for example, buffer 
overflow
probabilities, it would be useful to be able to ask questions about 
the probability of a given queue size ($\Prob{Q=i}$) or the
probability that the queue is more than a given size ($\Prob{Q > i}$).

The probability that the queue is zero is given by
\begin{equation*}
\Prob{Q=0} = \lim_{z \rightarrow 0} \E{z^Q},
\end{equation*}
and, more generally, the probability that the queue length is $q$ can
be found by differentiating $q$ times and taking the limit as 
$z \rightarrow 0$.
\begin{equation*}
\Prob{Q=i} = \lim_{z \rightarrow 0} \frac {d^q \E{z^Q}}{q! d^qz}.
\end{equation*}
This can be solved computationally by repeated
symbolic differentiation.  However, this is
computationally intensive and the algorithm is numerically unstable.
Another approach is to produce a recursive formula for the coefficient
of $z^i$ in $\E{z^Q}$.  This can be done using standard techniques
for formal power series from, for example, Knuth \cite{knuth1997}.

By definition, the coefficient of $z^i$ in $\E{z^Q}$ is
$\Prob{Q=i}$.  Let
$N_i$ and $D_i$ be the coefficients of $z^i$ in $N(z)$ and $D(z)$.
Since \eqref{eqn:ezq} is true 
for all $z$, therefore standard techniques for division of power 
series give
\begin{equation*}
\sum_{i=0}^k \Prob{Q = i} D_{k-i} = b_0 N_k,
\end{equation*}
which rearranges to
\begin{equation}
\Prob{Q = k} = \frac{1}{D_0} \left[N_kb_0 - \sum_{i=0}^{k-1} 
\Prob{Q=i} D_{k-i}\right].
\label{eqn:eqrec}
\end{equation}
This recursive formula expresses $\Prob{Q=k}$ in terms of $b_0$ which 
can
be evaluated with \eqref{eqn:b0}, coefficients $N_i$ and $D_i$ and 
$\Prob{Q=j}$
for $j < k$.  Therefore, the $\Prob{Q=k}$ can be calculated in turn 
beginning
with $\Prob{Q=0}$ which is given by
\begin{equation*}
\Prob{Q= 0} = \frac{b_0 N_0}{D_0}.
\end{equation*}

The coefficients $N_i$ and $D_i$ can be
easily calculated.  Let $G_{i,j}$ be the coefficient of $z^i$ in 
$(g(z)/z)^j$.  The
coefficients above can be expressed as
\begin{equation}
D_i= \delta_{i-1} -  \sum_{j=0}^n f_j G_{i,j}
\label{eqn:di}
\end{equation}
where $\delta_i$ is the Kronecker delta function 
($\delta_{i-1} = 1$ if $i=1$ and $0$ otherwise) and
\begin{equation}
N_i = 
\begin{cases}
- \sum_{j=0}^n G_{i,j} \sum_{k=j}^n f_k & i = 0 \\
\sum_{j=0}^n \left[G_{i-1,j} - G_{i,j}\right] \sum_{k=j}^n f_k & i > 
0.
\end{cases}
\label{eqn:ni}
\end{equation}
The coefficients $G_{i,j}$ can be calculated by another recurrence
relation. Since
\begin{equation*}
g(z)/z = \sum_{i=0}^{m-1} g_{i+1} z^i,
\end{equation*} 
and as $G_{i,0} = \delta_i$, this gives the recurrence relation
\begin{equation*}
G_{i,j+1} = \sum_{k=0}^i G_{k,j} g_{i+1-k}.
\end{equation*}

\subsection{A simpler model --- the constant batch size model}

A simplification occurs when the batch size is fixed.  Assume
that, work must arrive in units of exactly $r$ where $r > 1$ 
($r=1$ is the uninteresting system where no queue ever forms).
This means that $g(z)/z = z^{r-1}$ and $(g(z)/z)^j = z^{j(r-1)}$.
In turn this gives 
$G_{i,j} = \delta_{i-j(r-1)}$ where here, and throughout this
section, $\delta$ is the Kronecker delta function.
Obviously $\overline{g} = r$ and $\overline{g^2} = \overline{g}^2 = r^2$,
hence from \eqref{eqn:eqfinal}
$$
\E{Q} = \frac{r(r-1)(\overline{f^2} - \overline{f}^2)}
{2 ( 1 + \overline{f}) (1 + \overline{f} - r \overline{f})}.
$$
It can also be shown that
$
\Prob{Q=0} = b_0/f_0
$
and for $k > 0$,
\begin{align*}
\Prob{Q=k} & = \frac{1}{f_0}\biggl[ \Prob{Q=k-1} 
 - \sum_{j=0}^{k-1} \delta_{j/(r-1) - \lfloor j/(r-1) \rfloor} f_j \Prob{Q=j} \\
& - b_0 \delta_{(k-1)/(r-1) - \lfloor (k-1)/(r-1) \rfloor} \sum_{j=(k-1)/(r-1)}^n f_j \\
& + b_0 \delta_{k/(r-1) - \lfloor k/(r-1) \rfloor} \sum_{j=k/(r-1)}^n f_j 
\biggr],
\end{align*}
where $\lfloor x \rfloor$ is the floor function and
with the notational conveniences that $f_j= 0$ for
$j > n$ and that $\sum_{j}^n f_j = 0$ for $j > n$.  
Note that the delta functions here are simply
testing if a given expression, for example $k/(r-1)$, is integer.
From \eqref{eqn:b0},
$$
b_0= \frac{1+\overline{f} - r \overline{f}}{1 + \overline{f}}.
$$

\section{Simulation tests}
\label{sec:simulation}

The system needs to be tested against simulation to see if it can be
practically used.  While the equations from the previous sections are
correct, they are not useful if the numerical stability of the recursive
system of equations is poor.  For the expected queue length calculations
this is not an issue but it is for the probability of higher queue lengths
since it is likely that $\Prob{Q=i}$ will become extremely
small as $i$ becomes large.  This, in turn, will make the calculation
of \eqref{eqn:eqrec} problematic as $i$ becomes large.  Calculating
the potential effects of errors in the system of equations
given by \eqref{eqn:eqrec}, \eqref{eqn:di} and \eqref{eqn:ni} 
is non-trivial.  The answers presented here are tested
against simulation methods and appear valid for smaller queue
sizes but become obviously incorrect (negative probabilities for
example) for larger queue sizes.  The simulations here
were done in python.  The same calculations have
been tried using arbitrary precision arithmetic libraries.  
This enabled slightly larger queue sizes to be calculated
and give reasonable answers but at the expense of greatly increased
run time.

The method is simply to replicate the Markov chain and
queuing system described earlier and simulate it.
This is an exact simulation of the queuing system described
in the paper.
The simulation can then run for a set number of iterations and the 
queue
measured at each point to sample $\E{Q}$ or get a sample of the 
probabilities
$\Prob{Q=i}$ by measuring the proportion of the iterations where the 
queue
has the value $i$ in simulation.  Tables \ref{tab:parm1} or 
\ref{tab:parm2} show
the parameter sets for two different simulation
scenarios with the first representing a 
lightly loaded system and the second representing much heavier loading.
In both scenarios $\E{Q}$ matches well between theory and experiment as
expected and no results are presented here.

\begin{table}[htb]
\begin{center}
\begin{tabular}{|c |cccc|} \hline
$i$ & 0 & 1 & 2 & 3 \\ \hline
$f_i$ & 0.8 & 0.1 & 0.05 & 0.05 \\
$g_i$ & 0 & 0.4 &  0.4 & 0.2 \\ \hline
\end{tabular}
\end{center}
\caption{Parameter set 1 for simulation.}
\label{tab:parm1}
\end{table}

\begin{table}[htb]
\begin{center}
\begin{tabular}{|c |ccccc|} \hline
$i$ & 0 & 1 & 2 & 3 & 4\\ \hline
$f_i$ & 0.6 & 0.2 & 0.1 & 0.05 & 0.05 \\
$g_i$ & 0 & 0.2 & 0.6 & 0.1 & 0.1 \\ \hline
\end{tabular}
\end{center}
\caption{Parameter set 2 for simulation.}
\label{tab:parm2}
\end{table}

\begin{figure}[htb]
\begin{center}
\includegraphics[width=12cm]{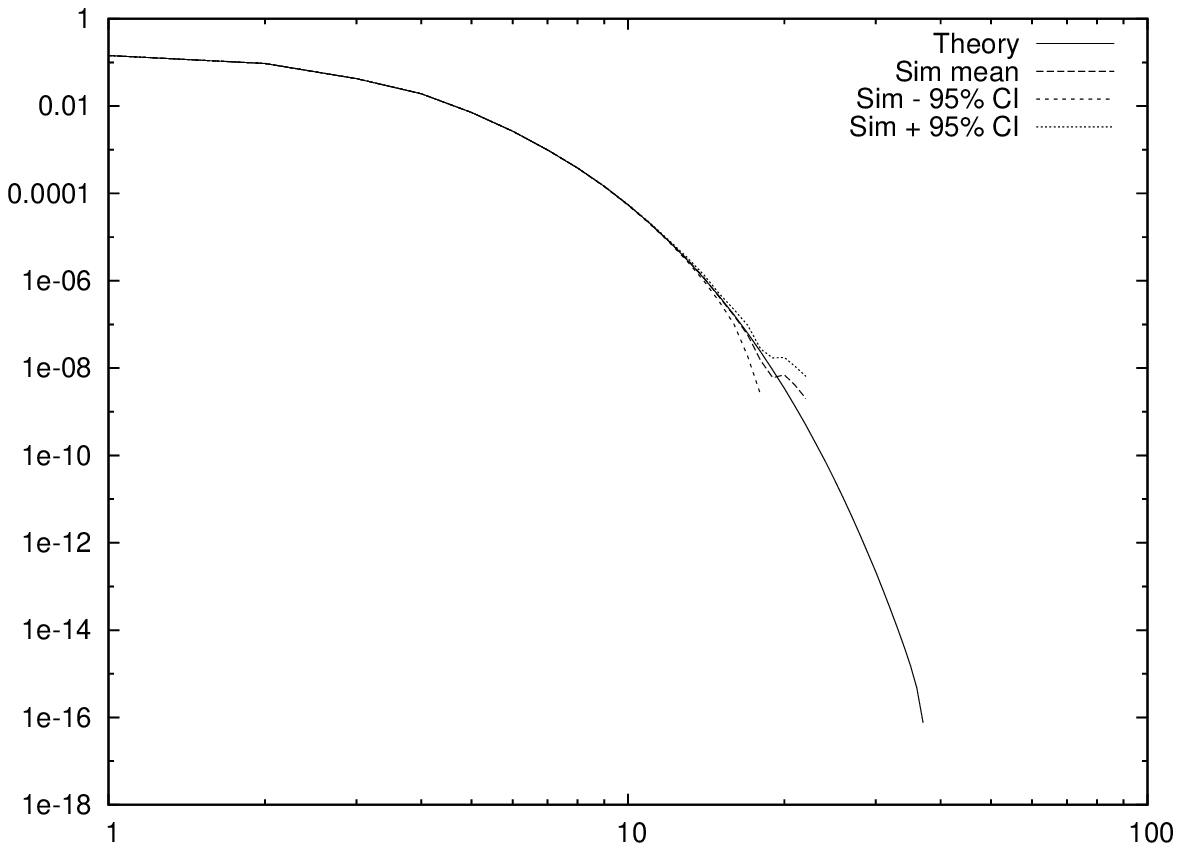}
\end{center}
\caption{Theory versus simulation for parameter set 1.}
\label{fig:parms1}
\end{figure}

\begin{figure}[htb]
\begin{center}
\includegraphics[width=12cm]{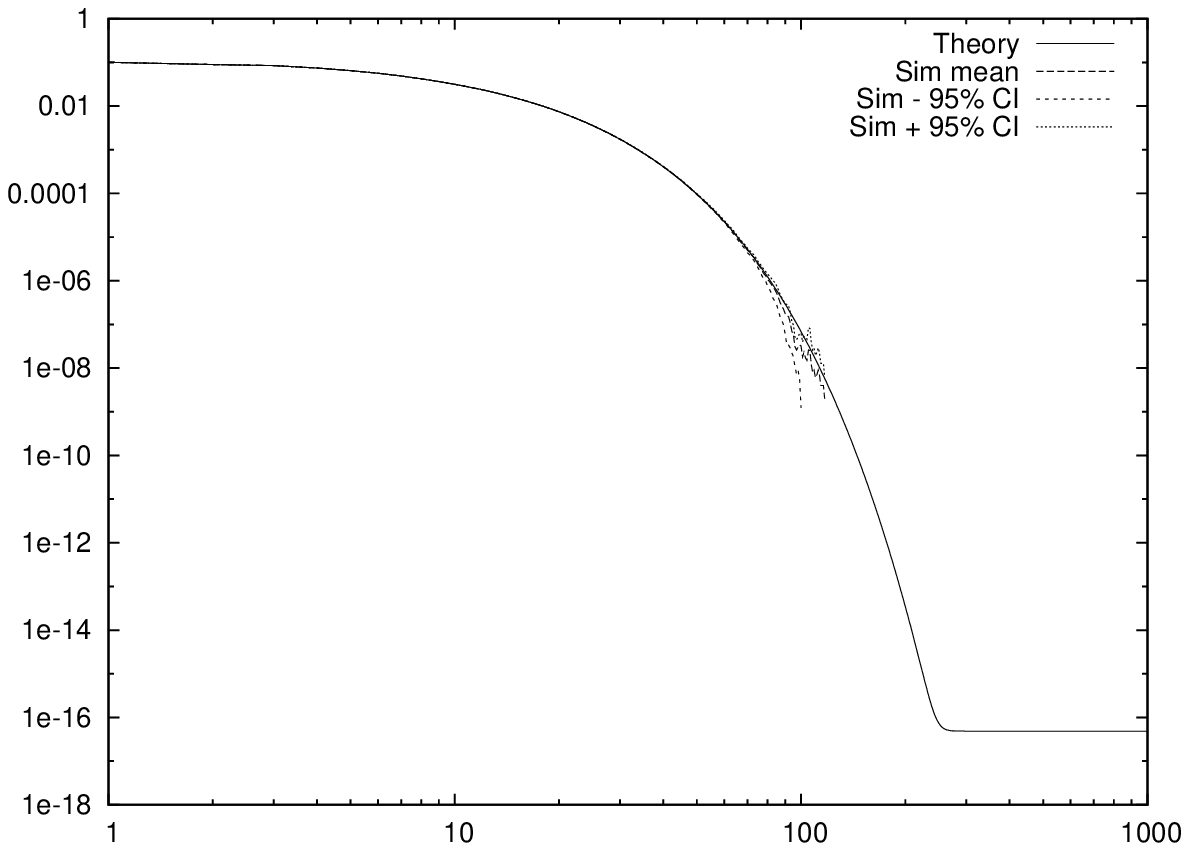}
\end{center}
\caption{Theory versus simulation for parameter set 2.}
\label{fig:parms2}
\end{figure}

Figures \ref{fig:parms1} and \ref{fig:parms2} show the theory from 
section
\ref{sec:pq} plotted against ten simulation runs for each parameter 
set, each run being $10^8$ iterations.  The simulation results
are shown as a mean for the ten simulation runs and upper and lower
95\% confidence intervals.  The plots are on a logscale and therefore, 
obviously, values which are zero or negative do not show up.  This
happens with lower confidence intervals and when there are
rounding errors in the theoretical calculation.

For
parameter set one, rounding errors make the theoretical
calculation obviously 
unreliable 
(this is obvious when the numbers become negative) at around
$i=40$ (where $\Prob{Q=i}$ is around $1\times 10^{-16}$) although it 
is likely that
some figures before this are rendered inaccurate due to rounding --- 
the negative
numbers are omitted from the plot because of the logscale, all others are included. 
In fact 
Note that because
only $10^8$ iterations were simulated the lowest sample 
probability that
can be assigned in simulation is $1\times 10^{-8}$.
For parameter set two, it can be seen from figure \ref{fig:parms2}
that the calculation remains reliable until above $i=110$ when,
again, it becomes obviously inaccurate due to rounding at
a probability again around $1 \times 10^{-16}$.
As
has been mentioned, by using arbitrary precision arithmetic
libraries, the numerical stability can be increased slightly 
but at the expense of greatly increased run times.

\section{Conclusions}

The arrival system described is quite general and could be useful in 
any system when work arrives at discrete times in discrete batches. 
The solutions given provide mean queue lengths and delays for the 
given arrival process.  In addition a system of equations has been given 
which can calculate the probability that the queue has a given length from the 
system parameters and the probabilities of smaller lengths.  
The numerical stability of the recursive system of equations 
giving the probability distribution has been assessed via simulation.

\bibliography{rgc_pe06}

\begin{thebibliography}{10}

\bibitem{barenco2004}
M.~Barenco and D.K. Arrowsmith.
\newblock The autocorrelation of double intermittency maps and the simulation
  of computer packet traffic.
\newblock {\em Dynamical Systems}, 19(1):61--74, 2004.

\bibitem{bernstein2000}
D.~S. Bernstein and C.~F.~Van Loan.
\newblock Rational matrix functions and rank-1 updates.
\newblock {\em SIAM J. Matrix Anal. Appl.}, 22(1):145--154, 2000.

\bibitem{blondia1992}
C.~Blondia and O.~Casals.
\newblock Statistical multiplexing of {VBR} sources: a matrix- analytic
  approach.
\newblock {\em Performance evaluation}, 16(1--3):5--20, 1992.

\bibitem{breuer2005}
L.~Breuer and D.~Baum.
\newblock {\em An introduction to queueing theory and matrix-analytic methods}.
\newblock Springer, 2005.

\bibitem{clegg2005}
R.~G. Clegg and M.~M. Dodson.
\newblock A {M}arkov based method for generating long-range dependence.
\newblock {\em Phys. Rev. E}, 72:026118, 2005.
\newblock Available online at: \\ {\tt
  www.richardclegg.org/pubs/rgcpre2004.pdf}.

\bibitem{coelho2000}
Z.~Coelho.
\newblock {\em Asymptotic laws for symbolic dynamical systems}, pages 123--165.
\newblock Number 279 in LMS Lecture Notes Series. Cambridge University Press,
  2000.

\bibitem{hermann1993}
Christoph Herrmann.
\newblock Analyis of the discrete time {SMP/D/1/s} finite buffer queue with
  applications in {ATM}.
\newblock {\em Proc. {IEEE} {INFOCOM}}, 1993.

\bibitem{khinchin1932}
A.~Y. Khinchin.
\newblock Mathematical theory of stationary queues.
\newblock {\em Mat. Sbornik}, 39:73--84, 1932.

\bibitem{klemm2003}
Alexander Klemm, Christoph Lindemann, and Marco Lohmann.
\newblock Modeling {IP} traffic using the batch {M}arkovian arrival process.
\newblock {\em Performance Evaluation}, 54:149--173, 2003.

\bibitem{knuth1997}
D.~E. Knuth.
\newblock {\em Seminumerical Algorithms}, volume~2 of {\em The Art of Computer
  Programming}, section 4.7, pages 525--533.
\newblock Addison-Wesley, Reading, Massachusetts, third edition, 1997.

\bibitem{fretwell1996}
Demetres~D. Kouvatsos and Rod~J. Fretwell.
\newblock Batch renewal process: exact model of traffic correlation.
\newblock In {\em Proceedings of the 2nd International Workshop on Architecture
  and Protocols for High Performance Networks}, pages 285--304, Deventer, The
  Netherlands, The Netherlands, 1996. Kluwer, B.V.

\bibitem{li1993}
S.~Li and C.~Hwang.
\newblock Queue response to input correlation functions: Discrete spectral
  analysis.
\newblock {\em {IEEE}\slash{ACM} Trans. on Networking}, 1(5):522--533, October
  1993.

\bibitem{little1961}
J.~Little.
\newblock A proof of the queueing formula {$L=\lambda W$}.
\newblock {\em Oper. Res. J.}, 18:172--174, 1961.

\bibitem{lucantoni1991}
D.~M. Lucantoni.
\newblock New results on the single server queue with a batch markovian arrival
  process.
\newblock {\em Stochastic models}, 7, 1991.

\bibitem{wang1989}
X.~J. Wang.
\newblock Statistical physics of temporal intermittency.
\newblock {\em Phys. Rev. A}, 40(11):6647--6661, 1989.

\bibitem{woolf2002}
M.~Woolf, D.~K. Arrowsmith, R.~J. Mongrag\'{o}n, and J.~M. Pitts.
\newblock Optimization and phase transitions in a chaotic model of data
  traffic.
\newblock {\em Phys. Rev. E}, 66:046106, 2002.

\end{thebibliography}

\end{document}